\documentclass[letterpaper, 10 pt, conference]{ieeeconf}  

\IEEEoverridecommandlockouts                              

\usepackage{cite}
\usepackage{amsmath,amssymb,amsfonts}
\usepackage{amsmath,bm}
\usepackage{algorithmic}
\usepackage{graphicx}
\usepackage{textcomp}
\usepackage{varioref}
\usepackage{wrapfig}
\usepackage{threeparttable}
\usepackage{dcolumn}
\newcolumntype{d}{D{.}{.}{-1}}
\makeglossary
\usepackage{fancyvrb}
\fvset{fontsize=\footnotesize,xleftmargin=2em}
\usepackage{lettrine}
\usepackage[colorlinks]{hyperref}
\usepackage{graphicx}
\usepackage{epstopdf}
\usepackage{attachfile}
\usepackage{adjustbox} 
\usepackage{anyfontsize}
\usepackage{booktabs}
\usepackage{array}
\usepackage{varwidth}
\usepackage[section]{placeins}
\usepackage[intoc]{nomencl}
\usepackage{ifthen}
\usepackage{relsize}
\usepackage{cleveref}
\graphicspath{{./Figures/}} %

\title{\LARGE \bf Entry Trajectory Optimization for Mars Science Laboratory Class Missions Using Indirect Unified Trigonometrization Method*
}

\author{Kshitij Mall$^{1}$ and Ehsan Taheri$^{1}$
\thanks{*This work was not supported by any organization.}
\thanks{$^{1}$Kshitij Mall and Ehsan Taheri are with the Department of Aerospace Engineering, Auburn University,
        Auburn, AL 36849, USA
        {\tt\small \{mall,etaheri\}@auburn.edu}}%
}

\begin{document}
\maketitle
\thispagestyle{empty}
\pagestyle{empty}
\graphicspath{{./Figures/}}
\bstctlcite{IEEEexample:BSTcontrol}
\begin{abstract}
Application of traditional indirect optimization methods to optimal control problems (OCPs) with control and state path constraints is not a straightforward task. However, recent advances in regularization techniques and numerical continuation methods have enabled application of indirect methods to very complex OCPs. This study demonstrates the utility and application of an advanced indirect method, the Unified Trigonometrization Method (UTM), to a Mars Science Laboratory type entry problem. The objective is  to maximize the parachute deployment altitude for a free-time, fixed-final-velocity entry trajectory. For entry vehicles, in addition to the bank angle that is characterized by bang-bang control profiles, there are typically three state path constraints that have to be considered, namely, the dynamic pressure, heat rate and g-load. This study shows that the UTM enables simultaneous regularization of the bang-bang control and satisfaction of the state path constraints. Two scenarios with and without state path constraints are considered. The results obtained using the UTM for both of these cases are found to be in excellent agreement with a direct optimization method. Furthermore, an interesting feature emerges in the optimal control profile of the UTM during the initial high-altitude part of the resulting optimal trajectory for the scenario with state path constraints, which has an appealing practical implication. 
\end{abstract}

\section{INTRODUCTION}
Planetary scientists conduct Mars exploration to find answers to important questions concerning life on Mars. The science experiments have to be transported to Mars and eventually landed on its surface through entry vehicles. One such entry mission was the Mars Science Laboratory (MSL), a robotic space mission that successfully landed \textit{Curiosity} rover on the surface of Mars in 2012 \cite{mendeck2014entry}. Entry, descent, and landing (EDL) phases pose difficult challenges to Mars missions. For an MSL-type mission, the maximum final  altitude is desired because it is typically a concern and a limitation to the use of parachutes or other aerodynamically-driven deceleration mechanisms \cite{garcia2008apollo,lafleur2011mars,zheng2017indirect}. For a preliminary design phase of a Mars mission, and before the selection of the landing site, the entry trajectory is often optimized and designed to achieve the highest parachute deployment altitude \cite{zheng2017indirect,lafleur2011mars,mendeck2002guidance}.

Atmospheric trajectory optimization is a challenging task \cite{li2014review,kamyar2014aircraft,berning2018rapid,kumar2018atmospheric}. Several direct and evolutionary optimization methods have been used to solve entry-type problems including particle swarm optimization \cite{grant2007mars}, genetic algorithm \cite{sorgenfrei2013exploration}, and a direct collocation method \cite{li2011mars}. Direct methods are frequently favored over the indirect methods due to 1) ease of implementation, 2) incorporation of different types of constraints, and 3) larger domain of convergence. However, the direct optimization methods trade optimality for convergence and may take significant computational effort to solve a complex OCP such as the Mars EDL trajectory optimization \cite{von1992direct,zheng2017indirect,grant2012rapid}. On the other hand, indirect optimization methods derive first-order necessary conditions of optimality that convert the original OCP into a boundary-value problem (BVP). The solution to the resulting BVPs are, at least, guaranteed to be local extremal. However, solution methods to these BVPs exhibit significant sensitivity to the unknowns (e.g., Lagrange multipliers and costates). Moreover, incorporation of state path constraints has been a constraining factor to the application of indirect methods to overly constrained OCPs. Nevertheless, direct and indirect optimization methods have their own advantageous and being able to use them for solving OCPs will provide additional insights as is demonstrated herein. 

The focus of this paper is not to draw a comparison between direct and indirect methods of optimization, but to present an alternative way to the existing frameworks based on indirect methods (such as the one proposed in Ref. \cite{zheng2017indirect}) for solving complex MSL-type entry trajectory optimization problems. The indirect framework used in this study is the Unified Trigonometrization Method (UTM), which has recently been used to solve complex atmospheric flight mechanics problems as shown in Ref. \cite{mall2020utmafm}. 

The two main contributions of this paper are: 1) to solve a complex Mars entry trajectory optimization problem with control and multiple state path constraints using the UTM (and with MATLAB's built-in solver \textit{bvp4c}), and 2) to showcase that the optimal control \textit{may} become non-unique in the early, high-altitude part of such an entry trajectory where the Martian atmosphere is too thin. 

The remainder of this paper is organized as follows. Section \ref{sec:marsedlutm} presents formulation of the two point boundary value problem (TPBVP) associated with the Mars entry trajectory optimization problem and the solution process. Section \ref{sec:results} demonstrates the results for two scenarios and draws a comparison with the solutions obtained using a direct pseudo-spectral method (PSM) \cite{patterson2014gpops}. Section \ref{sec:conc} includes the concluding remarks. 
\section{BOUNDARY-VALUE PROBLEM FORMULATION AND SOLUTION PROCESS}
\label{sec:marsedlutm}
\begin{figure}[!htbp]
\centering
\includegraphics[width=3.2in]{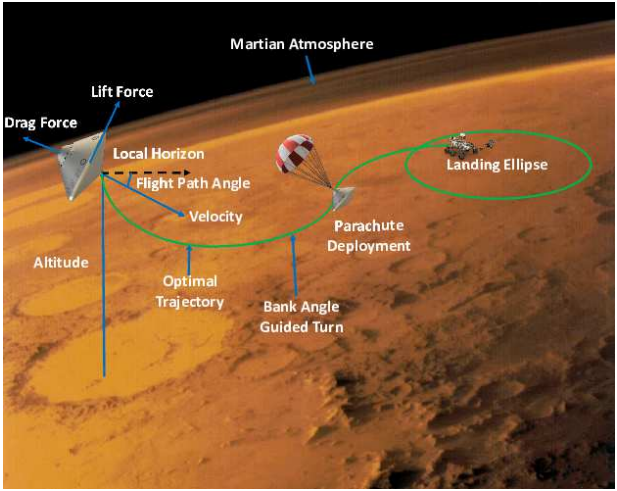}
\caption{Schematic of an MSL-type EDL on the Martian surface.}
\label{fig:msledlschem}
\end{figure}

We are considering a Mars entry problem with the objective of maximizing the final altitude at a predetermined final velocity from which the deceleration phase using a supersonic parachute is planned to be initiated. A schematic of an MSL-type EDL is shown in Fig.~\ref{fig:msledlschem} in which the objective is to land a rover on the surface of Mars in the vicinity of a designated landing site. This study, however, only concerns with the entry phase of the entire trajectory, which is the phase before the parachute deployment.  

The cost functional, $J$, is given (in Mayer form) in Eq.~\eqref{eqn:msledlobj}, while the equations of motion (EOMs) for a planar dynamics are given in Eqs.~(\ref{eqn:msledlobjeom1}--\ref{eqn:msledlobjeom3}). The resulting OCP assuming a non-rotational spherical Mars model with an exponential atmosphere and a point mass model for the entry vehicle can be written as
\begin{subequations}
\label{eqn:wec}
    \begin{align}
  J &= -h(t_{f}),\label{eqn:msledlobj}\\
  \dot{h} &= v\sin\gamma,\label{eqn:msledlobjeom1}\\
  \dot{v} &= -\dfrac{D}{m} - \dfrac{\mu\sin\gamma}{r^{2}},\label{eqn:msledlobjeom2}\\
  \dot{\gamma} &= \dfrac{L\cos\sigma}{mv} + \left(\dfrac{v}{r} - \dfrac{\mu}{r^{2}v}\right)\cos\gamma,\label{eqn:msledlobjeom3}
\end{align}
  \end{subequations}
where $r = r\textsubscript{MARS} + h$, $D = \dfrac{1}{2}\rho v^{2}C\textsubscript{D}A$, $L = \dfrac{1}{2}\rho v^{2}C\textsubscript{L}A$ and $\rho = \rho\textsubscript{0}\text{exp}(-h/H_{S})$. In these relations, \textit{r} is the radial magnitude, \textit{h} is the altitude, \textit{v} is the velocity magnitude, \textit{$\gamma$} is the flight path angle, \textit{m} is the mass of the entry vehicle, \textit{A} is the reference area of the vehicle, \textit{$\mu$} is the gravitational parameter of Mars, $r$\textsubscript{MARS} is Mars' mean radius, $\rho$\textsubscript{0} is the surface atmospheric density of Mars, \textit{D} is the drag force magnitude, \textit{L} is the lift force magnitude, \textit{$\sigma$} is the bank angle, and $H_S$ is the scale height. 

Three important along-the-path state constraints are imposed on entry trajectories based on the literature. These constraints are dynamic pressure, $q$, stagnation point heat-rate, $\dot{Q}$ \cite{sutton1971general}, and g-load \cite{yong2008rapid2,jorris2009three,grant2012rapid}, which can be written as
\begin{align}
q & = \dfrac{1}{2}\rho v^{2}\leq q_{\text{MAX}}, \label{eqn:marscons_q}\\
\dot{Q} & = k\sqrt{\dfrac{\rho}{r_{n}}}v^{3} \leq \dot{Q}_{\text{MAX}}, \label{eqn:marscons_Qdot}\\
g_{\text{LOAD}} &= \dfrac{\sqrt{L^{2}+D^{2}}}{mg}\leq g_{\text{MAX}}, \label{eqn:marscons_g}
\end{align}
where $q_{\text{MAX}}$ is 10 kPa, $\dot{Q}_{\text{MAX}}$ is 70 W/cm\textsuperscript{2}, $g_{\text{MAX}}$ is 5 Earth-gs, $k$ is a heat-rate constant, $r_{n}$ is the nose radius of the entry vehicle, and $g$ is 9.81 m/s\textsuperscript{2} \cite{li2011mars,jiang2016mars}. $\dot{Q}_{\text{MAX}}$ is governed by the limitations of the thermal protection system of the entry vehicle. On the other hand, $q_{\text{MAX}}$ and $g_{\text{MAX}}$ are limited due to the structural strength of the entry vehicle. Note that these maximum values upon the path constraints are corresponding to a robotic class mission. For a human-class mission, these constraints are usually much more stringent.  

Let $a_q\triangleq q/q_{\text{MAX}}$, $a_{\dot{Q}}\triangleq \dot{Q}/\dot{Q}_{\text{MAX}}$ and $a_{g}\triangleq g_{\text{LOAD}}/g_{\text{MAX}}$. The UTM regularizes this OCP using two simple trigonometric modifications given as
\begin{align}
&\cos\sigma = c_{0} + c_{1}\sin u,\label{eqn:msledlutmu}\\
&J = -h(t_{f}) + \int_{t_{0}}^{t_{f}}\left[\epsilon_{c}\cos u + \epsilon_{q}\sec\left(\frac{\pi}{2} a_q\right)\right]\text{d}t \nonumber\\
&\ \ + \int_{t_{0}}^{t_{f}}\left[\epsilon_{\dot{Q}}\sec\left(\dfrac{\pi}{2} a_{\dot{Q}}\right) + \epsilon_{g}\sec\left(\dfrac{\pi}{2} a_g\right)\right]\text{d}t,\label{eqn:msledlobjnew}
\end{align}
where $c_{0} = (u_{\text{MAX}} + u_{\text{MIN}})/2$ and $c_{1} = (u_{\text{MAX}} - u_{\text{MIN}})/2$ with $u_{\text{MIN}}$ and $u_{\text{MAX}}$ as $\cos(120^{\circ})$ and $\cos(30^{\circ})$, respectively.


Here, the original control input, $\sigma$, is expressed in terms of a new control input, $u$, as is defined in Eq.~\eqref{eqn:msledlutmu}. In addition, the state path constraints are incorporated into the cost functional through the running cost (i.e., the integrand also known as the Lagrangian). In addition, an error-control term, $\epsilon_c \cos(u)$, is also introduced into the Lagrangian. The idea of introducing regularized control and error-control terms was originally proposed by Silva and Trelat \cite{silva2010smooth}, and it was later enhanced by Mall and Grant \cite{mall2017epsilon} utilizing trigonometric terms.  

In Eq.~\eqref{eqn:msledlobjnew}, $\epsilon_{c}$ is an error parameter used to influence the amount of smoothing of the bang-bang control structure. Additionally, $\epsilon_{q}$, $\epsilon_{\dot{Q}}$, and $\epsilon_{g}$ are the weighting factors (continuation parameters) multiplied by the penalty terms, $\sec\left(\dfrac{\pi}{2} a_q \right)$, $\sec\left(\dfrac{\pi}{2} a_{\dot{Q}}\right)$, and $\sec\left(\dfrac{\pi}{2} a_g\right)$, respectively, to impose the three path constraints. When the value of the path constraint, for instance, $a_q \rightarrow 1$, the argument inside the first secant term in the Lagrangian approaches $\pi/2$, which results in an infinite value for the Lagrangian. Since the objective is to minimize these penalty terms, the solver avoids the maximum values for the path constraints, thereby leading to solutions that are extremely close to the exact optimal solutions. The Hamiltonian associated with the regularized OCP can be written as
\begin{align} 
H &= \lambda_{h}\dot{h} + \lambda_{v}\dot{v} + \lambda_{\gamma}\dot{\gamma} + \epsilon_{c}\cos u + \epsilon_{q}\sec\left(\dfrac{\pi}{2} a_q\right) \nonumber\\
& \ \ + \epsilon_{\dot{Q}}\sec\left(\dfrac{\pi}{2} a_{\dot{Q}}\right) + \epsilon_{g}\sec\left(\dfrac{\pi}{2} a_g\right).\label{eqn:msledlutmham}
\end{align}

The costates dynamics obtained using the Euler-Lagrange equation are
\begin{subequations}
\label{eqn:msledlutmcostates}
\begin{align}
    \dot{\lambda}_{h} &= \lambda_{\gamma}\left[\dfrac{Lu}{mv H_S} + \dfrac{\cos\gamma}{r^{2}}\left(v - \dfrac{2\mu}{rv}\right)\right] \nonumber\\
   & \quad - \lambda_{v}\left(\dfrac{D}{m H_S} + \dfrac{2\mu\sin\gamma}{r^{2}}\right) \nonumber\\
    & \quad + \dfrac{\pi\epsilon_{q}}{2H_S} a_q\sec\left(\dfrac{\pi}{2}a_q\right)\tan\left(\dfrac{\pi}{2}a_q\right) \nonumber \\
    & \quad + \dfrac{\pi\epsilon_{\dot{Q}}}{4\rho H_S}a_{\dot{Q}}\sec\left(\dfrac{\pi}{2}a_{\dot{Q}}\right)\tan\left(\dfrac{\pi}{2}a_{\dot{Q}}\right) \nonumber \\    
    & \quad + \dfrac{\pi\epsilon_{g}}{2H_S} a_g\sec\left(\dfrac{\pi}{2} a_g\right)\tan\left(\dfrac{\pi}{2}a_g\right),\label{eqn:msledlutmcostate1}\\
\dot{\lambda}_{v} &= -\lambda_{h}\sin\gamma + \dfrac{2\lambda_{v}D}{mv} - \dfrac{\lambda_{\gamma}Lu}{mv^{2}} \nonumber\\
& \quad - \dfrac{\pi\epsilon_{q}}{v}a_q\sec\left(\dfrac{\pi}{2}a_q\right)\tan\left(\dfrac{\pi}{2}a_q\right)\nonumber \\
& \quad - \dfrac{3\pi\epsilon_{\dot{Q}}}{2v}a_{\dot{Q}}\sec\left(\dfrac{\pi}{2}a_{\dot{Q}}\right)\tan\left(\dfrac{\pi}{2}a_{\dot{Q}}\right)\nonumber\\ 
& \quad - \dfrac{\pi\epsilon_{g}}{v}a_g\sec\left(\dfrac{\pi}{2}a_g\right)\tan\left(\dfrac{\pi}{2}a_g\right),\label{eqn:msledlutmcostate2}\\
\dot{\lambda}_{\gamma} &= \left ( -\lambda_{h}v + \dfrac{\lambda_{v}\mu}{r^{2}} \right ) \cos\gamma+ \lambda_{\gamma}\sin\gamma\left(\dfrac{v}{r} - \dfrac{\mu}{r^{2}v}\right). \label{eqn:msledlutmcostate3}
\end{align}
\end{subequations}

The TPBVP requires boundary conditions upon the costates and the Hamiltonian. For a free-final time problem, the final value of the Hamiltonian, $H(t_f)$, is 0. Additionally, the transversality condition on the altitude results in $\lambda_{h}(t_f) = -1$. Since all of the states are fixed at the initial point, the costates corresponding to these states are free and constitute part of the unknown variables. At the final point, only the velocity is known, and therefore $\lambda_{v}(t_f)$ is free. The remaining costate, $\lambda_{\gamma}(t_f)$, is 0 since $\gamma$ is free at the final time. 

The switching function for this problem, $H_{1}$, is shown in Eq.~\eqref{eqn:msledlutmswitch}. The strong form of optimality (i.e., $\partial H/\partial u = 0$) is used to obtain the optimal control law as shown in Eq.~\eqref{eqn:msledlutmtrigcontrol}, which is dependent on $H_{1}$, $\epsilon_{c}$, and $c_{1}$. Even if the value of $H_{1}$ vanishes, the optimal control can be explicitly found from among these two control options using the PMP, which states that the extremal control minimizes the Hamiltonian along an extremal trajectory.  
\begin{subequations}
\label{eqn:msledlutmtrigcontrolswitch}
\begin{align}
H_{1} &= 
\dfrac{\lambda_{\gamma}L}{mv},\label{eqn:msledlutmswitch}\\
u^{*} &= 
\begin{cases}
\arctan\left(\dfrac{c_{1}H_{1}}{\epsilon_{c}}\right),\\\\
\arctan\left(\dfrac{c_{1}H_{1}}{\epsilon_{c}}\right) + \pi.\label{eqn:msledlutmtrigcontrol}\\
\end{cases}
\end{align}
\end{subequations}

To implement the UTM, we have developed a framework in MATLAB with inputs as the state dynamics, cost functional, and boundary conditions. In an automated manner, this framework 1) generates the associated costate EOMs and the optimal control law using the Symbolic Math Toolbox of MATLAB, 2) generates the costate and Hamiltonian boundary conditions, and 3) implements the PMP.


\section{Results}
\label{sec:results}
This section is divided into two sub-sections: 1) a case in which the control constraint is considered while the state path constraints are ignored, and 2) a case in which the three state path constraints (see Eqs.~\eqref{eqn:marscons_q}-\eqref{eqn:marscons_g}) are taken into account in addition to the control constraint. The results obtained using the UTM are compared with a pseudo-spectral method (PSM) \cite{patterson2014gpops} to verify the accuracy of the results. The values for the constants used for numerical simulations are shown in Table~\ref{table:msledlconsts}, which are taken from Ref. \cite{zheng2017indirect}. Table~\ref{table:msledlbv} summarizes the boundary values for the two considered cases. 

 \begin{table}[h!t]
 \begin{center}
  \caption{Constants for the Mars entry problem.}
 \begin{tabular}{ccc}
 \hline
 \hline
 \textbf{Parameter} & \textbf{Unit} & \textbf{Value}\\\hline
 $r$\textsubscript{MARS} & km & 3397\\
 $\mu$ & km\textsuperscript{3}/s\textsuperscript{2} & 42840\\
 $\rho\textsubscript{0}$ & kg/m\textsuperscript{3} & 0.0158\\
 $H_S$ & km & 9.354\\
 $m$ & kg & 3300\\
 $A$ & m\textsuperscript{2} & 15.9\\
 $C$\textsubscript{D} & - & 1.45\\
 $C$\textsubscript{L} & - & 0.348\\
 $r_n$ & m & 0.6\\
 $k$ & kg\textsuperscript{1/2}/m\textsuperscript{2} & 1.9027$\times$10\textsuperscript{-4}\\
 \hline
\end{tabular}
 \label{table:msledlconsts}
 \end{center}
\end{table}

\begin{table}[h!tbp]
 \caption{Boundary conditions for the Mars entry problem.}
 \begin{center}{
 \begin{tabular}{cccc}
 \hline
 \hline
\textbf{Parameter} & \textbf{Unit} & \textbf{Initial Value} & \textbf{Final Value}\\\hline
$t$ & (s) & 0 & Free\\
$h$ & (km) & 125 & Free\\
$v$ & (km/s) & 6 & 0.54\\
$\gamma$ & (deg) & -11.5 & Free\\
$\lambda_{h}$ & (nd) & Free & -1\\
$\lambda_{v}$ & (s) & Free & Free\\
$\lambda_{\gamma}$ & (m/deg) & Free & 0\\
\hline
\end{tabular}}
\label{table:msledlbv}
 \end{center}
\end{table}

\subsection{Case I: Without State Path Constraints}
A numerical continuation approach\cite{grant2014rapid,taheri2016enhanced,junkins2018exploration} with two continuation sets has been adopted for this case. Using this continuation approach, a simpler OCP is solved initially with the following values: an initial $h = 50$ km, $t_f = 10$ seconds, and a high value for $\epsilon_{c}$ (1 m/s). The initial conditions for $v$ and $\gamma$ are not changed (see Table~\ref{table:msledlbv}). Note that a random initial guess value of -0.1 is chosen for all the costates for the simpler OCP. In the first continuation set, this simpler initial solution serves as an initial guess for a subsequent complex problem comprising a lower terminal value of $v$ and a higher initial value of $h$. The first continuation set is completed when the terminal $v$ and initial $h$ specified in Table~\ref{table:msledlbv} are reached after a specified number of steps. The subsequent continuation set operates on reducing $\epsilon_{c}$ to a reasonably small value of $1.0 \times 10^{-6}$. 

\begin{figure}[!htbp]
\centering
\includegraphics[width=3in]{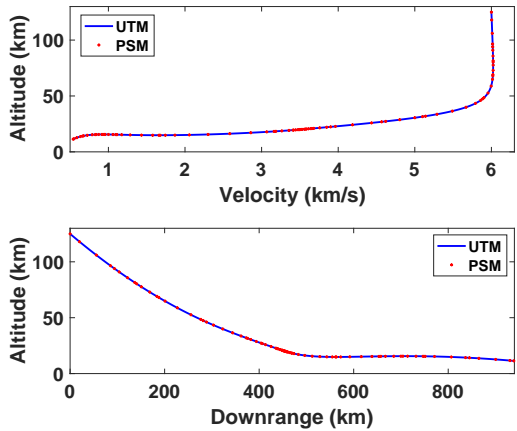}
\caption{Comparison of the energy (upper plot) and trajectory (lower plot) with the UTM and PSM for Case I.}
\label{fig:uc_energy_traj}
\end{figure}

Figure \ref{fig:uc_energy_traj} shows an excellent agreement between the results obtained using the UTM and PSM. The upper plot depicts the trajectory in the $(h$-$v)$-space (energy plot), where the vehicle trajectory initiates from the top-right point with the highest energy and terminates at the bottom-left point that corresponds to the least energy. The lower plot shows the change in altitude vs. downrange during the trajectory. The maximum terminal altitude attained for this case is $h_f = 11.3665$ km. For PSM $h_f = 11.3667$ km, which is extremely close to the result obtained using the UTM.

\begin{figure}[!htbp]
\centering
\includegraphics[width=3.5in]{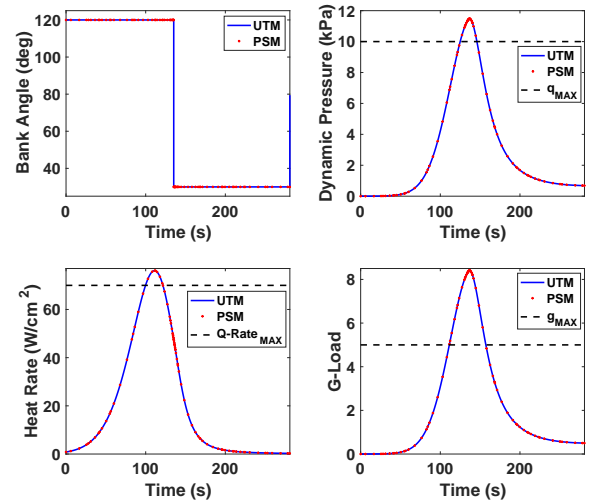}
\caption{Control and constraint plots with the UTM and PSM for Case I.}
\label{fig:uc_control_const}
\end{figure}

Figure \ref{fig:uc_control_const} contains the time histories of the values for the dynamic pressure, heat-rate, and g-load obtained using the UTM and PSM. The maximum allowable values for the state path constraints are denoted by dash lines. The results indicate that the constraints are clearly violated. The bank angle first stays at the maximum value to obtain a negative lift force making it possible for the vehicle to dive deeper into the atmosphere in order to gain higher magnitudes of the lift force. The vehicle reverses its bank angle and switches to its other extremum to have a positive lift force and consequently a higher terminal altitude. Note that the bank angle shifts from 30\textsuperscript{$\circ$} to 80\textsuperscript{$\circ$} at the very end. Since the bank angle at the terminal point does not impact the trajectory solution, this is an admissible result. 
\subsection{Case II: With Three State Path Constraints}
A numerical continuation approach similar to the Case I is used here with five continuation sets. The first continuation set is exactly similar to the previous case except that in this case we also need three state path constraint penalty parameters: $\epsilon_{Q}$, $\epsilon_{g}$, and $\epsilon_{q}$. Furthermore, the maximum values of the state path constraints are chosen initially as: $\dot{Q}_{\text{MAX}}$ = 200 W/cm\textsuperscript{2}, $g_{\text{MAX}}$ = 50 Earth-gs, and $q_{\text{MAX}}$ = 100 kPa, which result in none of the path constraints being active. A higher value of 1 m/s is chosen for $\epsilon_{c}$, $\epsilon_{Q}$, $\epsilon_{g}$, and $\epsilon_{q}$ while solving the first continuation set. The second continuation set runs upon $\dot{Q}_{\text{MAX}}$ and reduces its value to 70. The third and fourth continuation sets similarly bring down the values of $g_{\text{MAX}}$ and $q_{\text{MAX}}$ to 5 and 10, respectively. The fifth continuation step reduces the values of $\epsilon_{c}$, $\epsilon_{Q}$, $\epsilon_{g}$, and $\epsilon_{q}$ to $1.0 \times 10^{-6}$ each. 

\begin{figure}[!htbp]
\centering
\includegraphics[width=3in]{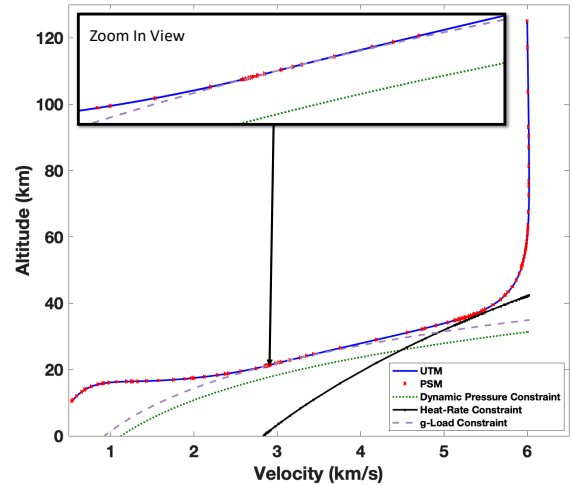}
\caption{Energy comparison plot for Case II.}
\label{fig:ac_energy}
\end{figure}

Figure \ref{fig:ac_energy} shows the energy plot obtained for this case using the UTM and PSM. The zoom in view shows that the g-load constraint becomes active, the heat-rate constraint becomes nearly active and the dynamic pressure constraint remains inactive throughout the entry phase. 

\begin{figure}[!htbp]
\centering
\includegraphics[width=3in]{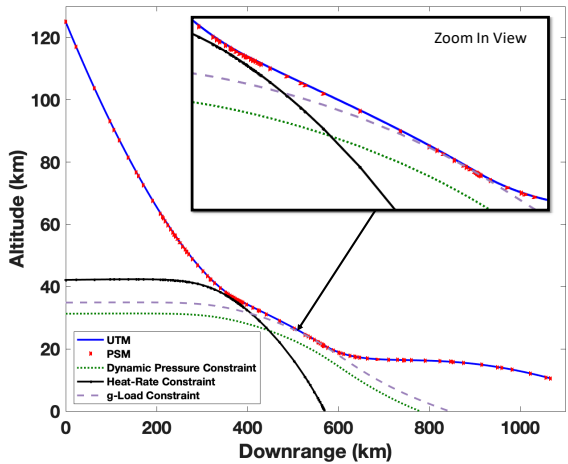}
\caption{Altitude vs. downrange comparison plot for Case II.}
\label{fig:ac_traj}
\end{figure}

The findings of the energy plot match with that of the trajectory plot shown in Fig.~\ref{fig:ac_traj}. The trajectory of the entry vehicle follows the g-load constraint path for a short time and nearly reaches the heat rate constraint at a certain instant as shown in the zoom in view. For this case, we obtain $h_f = 10.498$ km, which matches with the expectation and is slightly lesser than the case without any path constraints. This result matches exactly with that obtained using the PSM.

\begin{figure}[!htbp]
\centering
\includegraphics[width=3in]{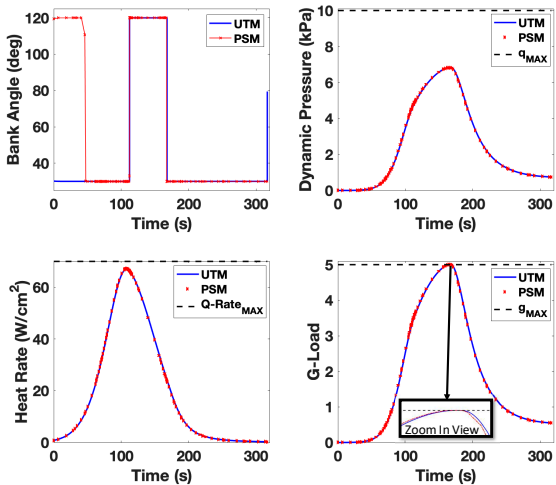}
\caption{Control and constraint comparison plots for Case II.}
\label{fig:ac_control_const}
\end{figure} 

Figure \ref{fig:ac_control_const} demonstrates that the UTM has been able to impose the control and the three state path constraints simultaneously. Furthermore, an interesting phenomenon is observed during the early parts of the control plot. The control solutions between the UTM and PSM match very well except for around first 47 seconds. Since the entry vehicle is within very thin regions of the atmosphere in the early part of the trajectory, the entry vehicle has negligible control forces. Thus, the control does not influence the dynamics of the problem and can become non-unique. Since there are only two controls possible in a bang-bang type control structure: the maximum or the minimum value, during the first 47 seconds of the trajectory, the control input can assume any of the possible values. 

\begin{figure}[!htbp]
\centering
\includegraphics[width=3.7in]{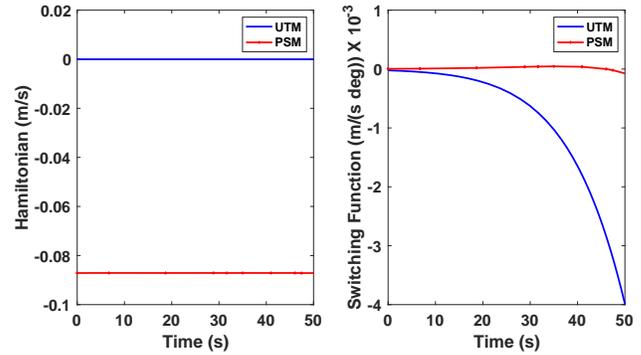}
\caption{First 50 seconds of the Hamiltonian and switching function time histories using the UTM and PSM for Case II.}
\label{fig:ham_sf_ac}
\end{figure}

This observation is further validated with the switching function plots shown in Fig.~\ref{fig:ham_sf_ac}. From an operational/practical point of view, the UTM solution is desirable as there are only two switches, which minimizes the control effort. Note that the bank angle again shifts at the last point while using the UTM similar to Case I. This is attributed to the fact that a regularized OCP with a minuscule error is solved using the UTM. The bank angle should in essence stay as 30\textsuperscript{$\circ$} at the final point of the trajectory without impacting the results. Therefore, the final switch in the bank angle can be ignored for the UTM.

The UTM results indicate that the entry vehicle, instead of diving down (using a negative lift force corresponding to the maximum bank angle value), should use the minimum bank angle corresponding to positive lift magnitudes. Thus, the entry vehicle avoids the denser atmosphere at higher velocities, thereby avoiding higher dynamic pressure, heat rate, and g-load paths, and following a safer trajectory while aiming at a maximum terminal altitude. As a consequence, the time of flight and downrange values both increase. Since the g-load path constraint is active, the entry vehicle attains a relatively lower altitude (7.6\% lower) compared to case I. 

\begin{figure}[!htbp]
\centering
\includegraphics[width=3in]{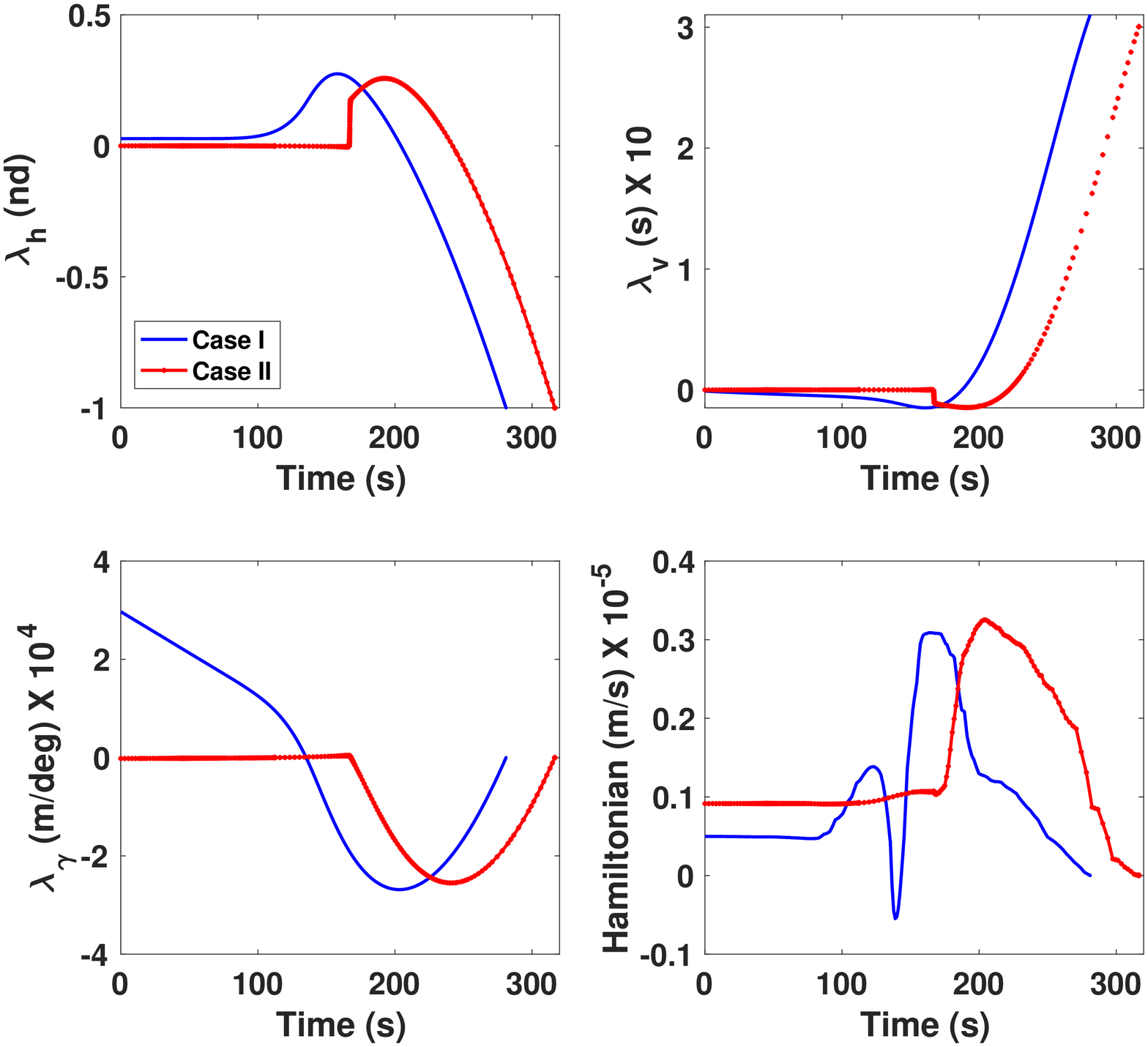}
\caption{Costates and Hamiltonian comparison plots for cases I and II obtained using the UTM.}
\label{fig:costates_utm_comp}
\end{figure}

The costates and the Hamiltonian obtained for the optimal trajectories of the two cases are compared in Fig.~\ref{fig:costates_utm_comp}. It is evident that since none of the path constraints are active for Case I, the costates are continuous. 

For Case II, the g-load constraint becomes active leading to a jump in $\lambda_{h}$ and $\lambda_{v}$ at the point where this constraint becomes active. The UTM is able to capture these jumps without solving a MPBVP as is done in the traditional indirect methods. The Hamiltonian for both the cases is continuous and close to 0 value throughout as shown in Fig.~\ref{fig:costates_utm_comp}.

The results for the two cases of the Mars entry trajectory optimization problem are summarized in Table~\ref{table:msledlresults}. Note that $q_{\text{PEAK}}$, $\dot{Q}_{\text{PEAK}}$, and $g$-$\text{load}_{\text{PEAK}}$ correspond to the peak values of the respective constraints on the optimal trajectories for both cases. All computations for were performed on a personal computer with a 2.6-GHz Intel i7 processor using MATLAB 2019a built-in BVP solver, \textit{bvp4c} (for the UTM). 
\begin{table}[!htbp]
          \caption{Comparison of the UTM results for the Mars entry problem.} 
          \centering
 \begin{tabular}{lcc}
 \hline
 \hline
        \textbf{Attribute} & \textbf{Case I} & \textbf{Case II}\\\hline
        \textbf{$\bm{h_{f}}$ (km)} & 11.367 & 10.498\\
        \textbf{Time of Flight (s)} & 280.999 & 316.607\\
        \textbf{Downrange (km)} & 938.813 & 1066.811\\
        \textbf{$\bm{\gamma_{f}}$ (deg)} & -13.083 & -13.996\\
        \textbf{$\bm{q_{\text{PEAK}}}$ (kPa)} & 11.478 & 6.825\\
        \textbf{$\bm{\dot{Q}_{\text{PEAK}}}$ (W/cm\textsuperscript{2})} & 76.123 & 67.028\\
        $\bm{g}$-$\bm{\textbf{\text{Load}}_{\textbf{\text{PEAK}}}}$ & 8.406 & 4.999\\\hline
\end{tabular}
          \label{table:msledlresults}
\end{table}

\section{CONCLUSIONS}
\label{sec:conc}

Application of a new indirect optimization framework, the Unified Trigonometrization Method (UTM), to two cases of a Mars Science Laboratory class entry problem were demonstrated. The objective was to maximize the final altitude from which a supersonic parachute will be deployed to further decelerate the vehicle. In Case I, a constraint was imposed only on the bank angle and the state path constraints (maximum dynamic pressure, heat rate, and g-load) were ignored. In Case II, the three aforementioned state path constraints were considered in addition to the constraint on the bank angle. The results were matched with those obtained using a pseudo-spectral direct method. High-quality results were obtained for the two cases by using the UTM. The final altitude obtained in Case II was found to be lesser than that obtained in Case I as expected. An interesting observation was made in Case II: \textit{the bank angle control in the initial (high-altitude) part of the constrained trajectory may become non-unique in case of Mars entry problems because the atmospheric forces are negligible}. This feature is evident especially during the first 47 seconds of the trajectory for Case II considered in this study. The control profile obtained using the UTM for Case II involves fewer bank angle reversals (control switches), which has practical utility and reduces the attitude control effort.

\bibliographystyle{IEEEtran} 
\bibliography{mslutm}

\end{document}